\documentclass[11pt]{article}
\long\def\ig#1{\relax}
\ig{Thanks to Roberto Minio for this def'n.  Compare the def'n of
\comment in AMSTeX.}

\newcount \coefa
\newcount \coefb
\newcount \coefc
\newcount\tempcounta
\newcount\tempcountb
\newcount\tempcountc
\newcount\tempcountd
\newcount\xext
\newcount\yext
\newcount\xoff
\newcount\yoff
\newcount\gap%
\newcount\arrowtypea
\newcount\arrowtypeb
\newcount\arrowtypec
\newcount\arrowtyped
\newcount\arrowtypee
\newcount\height
\newcount\width
\newcount\xpos
\newcount\ypos
\newcount\run
\newcount\rise
\newcount\arrowlength
\newcount\halflength
\newcount\arrowtype
\newdimen\tempdimen
\newdimen\xlen
\newdimen\ylen
\newsavebox{\tempboxa}%
\newsavebox{\tempboxb}%
\newsavebox{\tempboxc}%

\makeatletter
\def\settypes(#1,#2,#3){\arrowtypea#1 \arrowtypeb#2 \arrowtypec#3}
\def\settoheight#1#2{\setbox\@tempboxa\hbox{#2}#1\ht\@tempboxa\relax}%
\def\settodepth#1#2{\setbox\@tempboxa\hbox{#2}#1\dp\@tempboxa\relax}%
\def\settokens[#1`#2`#3`#4]{%
     \def\tokena{#1}\def\tokenb{#2}\def\tokenc{#3}\def\tokend{#4}}
\def\setsqparms[#1`#2`#3`#4;#5`#6]{%
\arrowtypea #1
\arrowtypeb #2
\arrowtypec #3
\arrowtyped #4
\width #5
\height #6
}
\def\setpos(#1,#2){\xpos=#1 \ypos#2}

\def\bfig{\begin{picture}(\xext,\yext)(\xoff,\yoff)}
\def\efig{\end{picture}}

\def\putbox(#1,#2)#3{\put(#1,#2){\makebox(0,0){$#3$}}}

\def\settriparms[#1`#2`#3;#4]{\settripairparms[#1`#2`#3`1`1;#4]}%

\def\settripairparms[#1`#2`#3`#4`#5;#6]{%
\arrowtypea #1
\arrowtypeb #2
\arrowtypec #3
\arrowtyped #4
\arrowtypee #5
\width #6
\height #6
}

\def\resetparms{\settripairparms[1`1`1`1`1;500]\width 500}

\resetparms

\def\mvector(#1,#2)#3{
\put(0,0){\vector(#1,#2){#3}}%
\put(0,0){\vector(#1,#2){30}}%
}
\def\evector(#1,#2)#3{{
\arrowlength #3
\put(0,0){\vector(#1,#2){\arrowlength}}%
\advance \arrowlength by-30
\put(0,0){\vector(#1,#2){\arrowlength}}%
}}

\def\horsize#1#2{%
\settowidth{\tempdimen}{$#2$}%
#1=\tempdimen
\divide #1 by\unitlength
}

\def\vertsize#1#2{%
\settoheight{\tempdimen}{$#2$}%
#1=\tempdimen
\settodepth{\tempdimen}{$#2$}%
\advance #1 by\tempdimen
\divide #1 by\unitlength
}

\def\vertadjust[#1`#2`#3]{%
\vertsize{\tempcounta}{#1}%
\vertsize{\tempcountb}{#2}%
\ifnum \tempcounta<\tempcountb \tempcounta=\tempcountb \fi
\divide\tempcounta by2
\vertsize{\tempcountb}{#3}%
\ifnum \tempcountb>0 \advance \tempcountb by20 \fi
\ifnum \tempcounta<\tempcountb \tempcounta=\tempcountb \fi
}

\def\horadjust[#1`#2`#3]{%
\horsize{\tempcounta}{#1}%
\horsize{\tempcountb}{#2}%
\ifnum \tempcounta<\tempcountb \tempcounta=\tempcountb \fi
\divide\tempcounta by20
\horsize{\tempcountb}{#3}%
\ifnum \tempcountb>0 \advance \tempcountb by60 \fi
\ifnum \tempcounta<\tempcountb \tempcounta=\tempcountb \fi
}

\ig{ In this procedure, #1 is the paramater that sticks out all the way,
#2 sticks out the least and #3 is a label sticking out half way.  #4 is
the amount of the offset.}

\def\sladjust[#1`#2`#3]#4{%
\tempcountc=#4
\horsize{\tempcounta}{#1}%
\divide \tempcounta by2
\horsize{\tempcountb}{#2}%
\divide \tempcountb by2
\advance \tempcountb by-\tempcountc
\ifnum \tempcounta<\tempcountb \tempcounta=\tempcountb\fi
\divide \tempcountc by2
\horsize{\tempcountb}{#3}%
\advance \tempcountb by-\tempcountc
\ifnum \tempcountb>0 \advance \tempcountb by80\fi
\ifnum \tempcounta<\tempcountb \tempcounta=\tempcountb\fi
\advance\tempcounta by20
}

\def\putvector(#1,#2)(#3,#4)#5#6{{%
\xpos=#1
\ypos=#2
\run=#3
\rise=#4
\arrowlength=#5
\arrowtype=#6
\ifnum \arrowtype<0
    \ifnum \run=0
        \advance \ypos by-\arrowlength
    \else
        \tempcounta \arrowlength
        \multiply \tempcounta by\rise
        \divide \tempcounta by\run
        \ifnum\run>0
            \advance \xpos by\arrowlength
            \advance \ypos by\tempcounta
        \else
            \advance \xpos by-\arrowlength
            \advance \ypos by-\tempcounta
        \fi
    \fi
    \multiply \arrowtype by-1
    \multiply \rise by-1
    \multiply \run by-1
\fi
\ifnum \arrowtype=1
    \put(\xpos,\ypos){\vector(\run,\rise){\arrowlength}}%
\else\ifnum \arrowtype=2
    \put(\xpos,\ypos){\mvector(\run,\rise)\arrowlength}%
\else\ifnum\arrowtype=3
    \put(\xpos,\ypos){\evector(\run,\rise){\arrowlength}}%
\fi\fi\fi
}}

\def\putsplitvector(#1,#2)#3#4{
\xpos #1
\ypos #2
\arrowtype #4
\halflength #3
\arrowlength #3
\gap 140
\advance \halflength by-\gap
\divide \halflength by2
\ifnum \arrowtype=1
    \put(\xpos,\ypos){\line(0,-1){\halflength}}%
    \advance\ypos by-\halflength
    \advance\ypos by-\gap
    \put(\xpos,\ypos){\vector(0,-1){\halflength}}%
\else\ifnum \arrowtype=2
    \put(\xpos,\ypos){\line(0,-1)\halflength}%
    \put(\xpos,\ypos){\vector(0,-1)3}%
    \advance\ypos by-\halflength
    \advance\ypos by-\gap
    \put(\xpos,\ypos){\vector(0,-1){\halflength}}%
\else\ifnum\arrowtype=3
    \put(\xpos,\ypos){\line(0,-1)\halflength}%
    \advance\ypos by-\halflength
    \advance\ypos by-\gap
    \put(\xpos,\ypos){\evector(0,-1){\halflength}}%
\else\ifnum \arrowtype=-1
    \advance \ypos by-\arrowlength
    \put(\xpos,\ypos){\line(0,1){\halflength}}%
    \advance\ypos by\halflength
    \advance\ypos by\gap
    \put(\xpos,\ypos){\vector(0,1){\halflength}}%
\else\ifnum \arrowtype=-2
    \advance \ypos by-\arrowlength
    \put(\xpos,\ypos){\line(0,1)\halflength}%
    \put(\xpos,\ypos){\vector(0,1)3}%
    \advance\ypos by\halflength
    \advance\ypos by\gap
    \put(\xpos,\ypos){\vector(0,1){\halflength}}%
\else\ifnum\arrowtype=-3
    \advance \ypos by-\arrowlength
    \put(\xpos,\ypos){\line(0,1)\halflength}%
    \advance\ypos by\halflength
    \advance\ypos by\gap
    \put(\xpos,\ypos){\evector(0,1){\halflength}}%
\fi\fi\fi\fi\fi\fi
}

\def\putmorphism(#1)(#2,#3)[#4`#5`#6]#7#8#9{{%
\run #2
\rise #3
\ifnum\rise=0
  \puthmorphism(#1)[#4`#5`#6]{#7}{#8}{#9}%
\else\ifnum\run=0
  \putvmorphism(#1)[#4`#5`#6]{#7}{#8}{#9}%
\else
\setpos(#1)%
\arrowlength #7
\arrowtype #8
\ifnum\run=0
\else\ifnum\rise=0
\else
\ifnum\run>0
    \coefa=1
\else
   \coefa=-1
\fi
\ifnum\arrowtype>0
   \coefb=0
   \coefc=-1
\else
   \coefb=\coefa
   \coefc=1
   \arrowtype=-\arrowtype
\fi
\width=2
\multiply \width by\run
\divide \width by\rise
\ifnum \width<0  \width=-\width\fi
\advance\width by60
\if l#9 \width=-\width\fi
\putbox(\xpos,\ypos){#4}
{\multiply \coefa by\arrowlength
\advance\xpos by\coefa
\multiply \coefa by\rise
\divide \coefa by\run
\advance \ypos by\coefa
\putbox(\xpos,\ypos){#5} }%
{\multiply \coefa by\arrowlength
\divide \coefa by2
\advance \xpos by\coefa
\advance \xpos by\width
\multiply \coefa by\rise
\divide \coefa by\run
\advance \ypos by\coefa
\if l#9%
   \put(\xpos,\ypos){\makebox(0,0)[r]{$#6$}}%
\else\if r#9%
   \put(\xpos,\ypos){\makebox(0,0)[l]{$#6$}}%
\fi\fi }%
{\multiply \rise by-\coefc
\multiply \run by-\coefc
\multiply \coefb by\arrowlength
\advance \xpos by\coefb
\multiply \coefb by\rise
\divide \coefb by\run
\advance \ypos by\coefb
\multiply \coefc by70
\advance \ypos by\coefc
\multiply \coefc by\run
\divide \coefc by\rise
\advance \xpos by\coefc
\multiply \coefa by140
\multiply \coefa by\run
\divide \coefa by\rise
\advance \arrowlength by\coefa
\ifnum \arrowtype=1
   \put(\xpos,\ypos){\vector(\run,\rise){\arrowlength}}%
\else\ifnum\arrowtype=2
   \put(\xpos,\ypos){\mvector(\run,\rise){\arrowlength}}%
\else\ifnum\arrowtype=3
   \put(\xpos,\ypos){\evector(\run,\rise){\arrowlength}}%
\fi\fi\fi}\fi\fi\fi\fi}}

\def\puthmorphism(#1,#2)[#3`#4`#5]#6#7#8{{%
\xpos #1
\ypos #2
\width #6
\arrowlength #6
\putbox(\xpos,\ypos){#3\vphantom{#4}}%
{\advance \xpos by\arrowlength
\putbox(\xpos,\ypos){\vphantom{#3}#4}}%
\horsize{\tempcounta}{#3}%
\horsize{\tempcountb}{#4}%
\divide \tempcounta by2
\divide \tempcountb by2
\advance \tempcounta by30
\advance \tempcountb by30
\advance \xpos by\tempcounta
\advance \arrowlength by-\tempcounta
\advance \arrowlength by-\tempcountb
\putvector(\xpos,\ypos)(1,0){\arrowlength}{#7}%
\divide \arrowlength by2
\advance \xpos by\arrowlength
\vertsize{\tempcounta}{#5}%
\divide\tempcounta by2
\advance \tempcounta by20
\if a#8 %
   \advance \ypos by\tempcounta
   \putbox(\xpos,\ypos){#5}%
\else
   \advance \ypos by-\tempcounta
   \putbox(\xpos,\ypos){#5}%
\fi}}

\def\putvmorphism(#1,#2)[#3`#4`#5]#6#7#8{{%
\xpos #1
\ypos #2
\arrowlength #6
\arrowtype #7
\settowidth{\xlen}{$#5$}%
\putbox(\xpos,\ypos){#3}%
{\advance \ypos by-\arrowlength
\putbox(\xpos,\ypos){#4}}%
{\advance\arrowlength by-140
\advance \ypos by-70
\ifdim\xlen>0pt
   \if m#8%
      \putsplitvector(\xpos,\ypos){\arrowlength}{\arrowtype}%
   \else
      \putvector(\xpos,\ypos)(0,-1){\arrowlength}{\arrowtype}%
   \fi
\else
   \putvector(\xpos,\ypos)(0,-1){\arrowlength}{\arrowtype}%
\fi}%
\ifdim\xlen>0pt
   \divide \arrowlength by2
   \advance\ypos by-\arrowlength
   \if l#8%
      \advance \xpos by-40
      \put(\xpos,\ypos){\makebox(0,0)[r]{$#5$}}%
   \else\if r#8%
      \advance \xpos by40
      \put(\xpos,\ypos){\makebox(0,0)[l]{$#5$}}%
   \else
      \putbox(\xpos,\ypos){#5}%
   \fi\fi
\fi
}}

\def\topadjust[#1`#2`#3]{%
\yoff=10
\vertadjust[#1`#2`{#3}]%
\advance \yext by\tempcounta
\advance \yext by 10
}
\def\botadjust[#1`#2`#3]{%
\vertadjust[#1`#2`{#3}]%
\advance \yext by\tempcounta
\advance \yoff by-\tempcounta
}
\def\leftadjust[#1`#2`#3]{%
\xoff=0
\horadjust[#1`#2`{#3}]%
\advance \xext by\tempcounta
\advance \xoff by-\tempcounta
}
\def\rightadjust[#1`#2`#3]{%
\horadjust[#1`#2`{#3}]%
\advance \xext by\tempcounta
}
\def\rightsladjust[#1`#2`#3]{%
\sladjust[#1`#2`{#3}]{\width}%
\advance \xext by\tempcounta
}
\def\leftsladjust[#1`#2`#3]{%
\xoff=0
\sladjust[#1`#2`{#3}]{\width}%
\advance \xext by\tempcounta
\advance \xoff by-\tempcounta
}
\def\adjust[#1`#2;#3`#4;#5`#6;#7`#8]{%
\topadjust[#1``{#2}]
\leftadjust[#3``{#4}]
\rightadjust[#5``{#6}]
\botadjust[#7``{#8}]}

\def\putsquarep<#1>(#2)[#3;#4`#5`#6`#7]{{%
\setsqparms[#1]%
\setpos(#2)%
\settokens[#3]%
\puthmorphism(\xpos,\ypos)[\tokenc`\tokend`{#7}]{\width}{\arrowtyped}b%
\advance\ypos by \height
\puthmorphism(\xpos,\ypos)[\tokena`\tokenb`{#4}]{\width}{\arrowtypea}a%
\putvmorphism(\xpos,\ypos)[``{#5}]{\height}{\arrowtypeb}l%
\advance\xpos by \width
\putvmorphism(\xpos,\ypos)[``{#6}]{\height}{\arrowtypec}r%
}}

\def\putsquare{\@ifnextchar <{\putsquarep}{\putsquarep%
   <\arrowtypea`\arrowtypeb`\arrowtypec`\arrowtyped;\width`\height>}}
\def\square{\@ifnextchar< {\squarep}{\squarep
   <\arrowtypea`\arrowtypeb`\arrowtypec`\arrowtyped;\width`\height>}}
\def\squarep<#1>[#2`#3`#4`#5;#6`#7`#8`#9]{{
\setsqparms[#1]
\xext=\width                                          
\yext=\height                                         
\topadjust[#2`#3`{#6}]
\botadjust[#4`#5`{#9}]
\leftadjust[#2`#4`{#7}]
\rightadjust[#3`#5`{#8}]
\begin{picture}(\xext,\yext)(\xoff,\yoff)
\putsquarep<\arrowtypea`\arrowtypeb`\arrowtypec`\arrowtyped;\width`\height>%
(0,0)[#2`#3`#4`#5;#6`#7`#8`{#9}]%
\end{picture}%
}}

\def\putptrianglep<#1>(#2,#3)[#4`#5`#6;#7`#8`#9]{{%
\settriparms[#1]%
\xpos=#2 \ypos=#3
\advance\ypos by \height
\puthmorphism(\xpos,\ypos)[#4`#5`{#7}]{\height}{\arrowtypea}a%
\putvmorphism(\xpos,\ypos)[`#6`{#8}]{\height}{\arrowtypeb}l%
\advance\xpos by\height
\putmorphism(\xpos,\ypos)(-1,-1)[``{#9}]{\height}{\arrowtypec}r%
}}

\def\putptriangle{\@ifnextchar <{\putptrianglep}{\putptrianglep
   <\arrowtypea`\arrowtypeb`\arrowtypec;\height>}}
\def\ptriangle{\@ifnextchar <{\ptrianglep}{\ptrianglep
   <\arrowtypea`\arrowtypeb`\arrowtypec;\height>}}

\def\ptrianglep<#1>[#2`#3`#4;#5`#6`#7]{{
\settriparms[#1]%
\width=\height                         
\xext=\width                           
\yext=\width                           
\topadjust[#2`#3`{#5}]
\botadjust[#3``]
\leftadjust[#2`#4`{#6}]
\rightsladjust[#3`#4`{#7}]
\begin{picture}(\xext,\yext)(\xoff,\yoff)
\putptrianglep<\arrowtypea`\arrowtypeb`\arrowtypec;\height>%
(0,0)[#2`#3`#4;#5`#6`{#7}]%
\end{picture}%
}}

\def\putqtrianglep<#1>(#2,#3)[#4`#5`#6;#7`#8`#9]{{%
\settriparms[#1]%
\xpos=#2 \ypos=#3
\advance\ypos by\height
\puthmorphism(\xpos,\ypos)[#4`#5`{#7}]{\height}{\arrowtypea}a%
\putmorphism(\xpos,\ypos)(1,-1)[``{#8}]{\height}{\arrowtypeb}l%
\advance\xpos by\height
\putvmorphism(\xpos,\ypos)[`#6`{#9}]{\height}{\arrowtypec}r%
}}

\def\putqtriangle{\@ifnextchar <{\putqtrianglep}{\putqtrianglep
   <\arrowtypea`\arrowtypeb`\arrowtypec;\height>}}
\def\qtriangle{\@ifnextchar <{\qtrianglep}{\qtrianglep
   <\arrowtypea`\arrowtypeb`\arrowtypec;\height>}}

\def\qtrianglep<#1>[#2`#3`#4;#5`#6`#7]{{
\settriparms[#1]
\width=\height                         
\xext=\width                           
\yext=\height                          
\topadjust[#2`#3`{#5}]
\botadjust[#4``]
\leftsladjust[#2`#4`{#6}]
\rightadjust[#3`#4`{#7}]
\begin{picture}(\xext,\yext)(\xoff,\yoff)
\putqtrianglep<\arrowtypea`\arrowtypeb`\arrowtypec;\height>%
(0,0)[#2`#3`#4;#5`#6`{#7}]%
\end{picture}%
}}

\def\putdtrianglep<#1>(#2,#3)[#4`#5`#6;#7`#8`#9]{{%
\settriparms[#1]%
\xpos=#2 \ypos=#3
\puthmorphism(\xpos,\ypos)[#5`#6`{#9}]{\height}{\arrowtypec}b%
\advance\xpos by \height \advance\ypos by\height
\putmorphism(\xpos,\ypos)(-1,-1)[``{#7}]{\height}{\arrowtypea}l%
\putvmorphism(\xpos,\ypos)[#4``{#8}]{\height}{\arrowtypeb}r%
}}

\def\putdtriangle{\@ifnextchar <{\putdtrianglep}{\putdtrianglep
   <\arrowtypea`\arrowtypeb`\arrowtypec;\height>}}
\def\dtriangle{\@ifnextchar <{\dtrianglep}{\dtrianglep
   <\arrowtypea`\arrowtypeb`\arrowtypec;\height>}}

\def\dtrianglep<#1>[#2`#3`#4;#5`#6`#7]{{
\settriparms[#1]
\width=\height                         
\xext=\width                           
\yext=\height                          
\topadjust[#2``]
\botadjust[#3`#4`{#7}]
\leftsladjust[#3`#2`{#5}]
\rightadjust[#2`#4`{#6}]
\begin{picture}(\xext,\yext)(\xoff,\yoff)
\putdtrianglep<\arrowtypea`\arrowtypeb`\arrowtypec;\height>%
(0,0)[#2`#3`#4;#5`#6`{#7}]%
\end{picture}%
}}

\def\putbtrianglep<#1>(#2,#3)[#4`#5`#6;#7`#8`#9]{{%
\settriparms[#1]%
\xpos=#2 \ypos=#3
\puthmorphism(\xpos,\ypos)[#5`#6`{#9}]{\height}{\arrowtypec}b%
\advance\ypos by\height
\putmorphism(\xpos,\ypos)(1,-1)[``{#8}]{\height}{\arrowtypeb}r%
\putvmorphism(\xpos,\ypos)[#4``{#7}]{\height}{\arrowtypea}l%
}}

\def\putbtriangle{\@ifnextchar <{\putbtrianglep}{\putbtrianglep
   <\arrowtypea`\arrowtypeb`\arrowtypec;\height>}}
\def\btriangle{\@ifnextchar <{\btrianglep}{\btrianglep
   <\arrowtypea`\arrowtypeb`\arrowtypec;\height>}}

\def\btrianglep<#1>[#2`#3`#4;#5`#6`#7]{{
\settriparms[#1]
\width=\height                         
\xext=\width                           
\yext=\height                          
\topadjust[#2``]
\botadjust[#3`#4`{#7}]
\leftadjust[#2`#3`{#5}]
\rightsladjust[#4`#2`{#6}]
\begin{picture}(\xext,\yext)(\xoff,\yoff)
\putbtrianglep<\arrowtypea`\arrowtypeb`\arrowtypec;\height>%
(0,0)[#2`#3`#4;#5`#6`{#7}]%
\end{picture}%
}}

\def\putAtrianglep<#1>(#2,#3)[#4`#5`#6;#7`#8`#9]{{%
\settriparms[#1]%
\xpos=#2 \ypos=#3
{\multiply \height by2
\puthmorphism(\xpos,\ypos)[#5`#6`{#9}]{\height}{\arrowtypec}b}%
\advance\xpos by\height \advance\ypos by\height
\putmorphism(\xpos,\ypos)(-1,-1)[#4``{#7}]{\height}{\arrowtypea}l%
\putmorphism(\xpos,\ypos)(1,-1)[``{#8}]{\height}{\arrowtypeb}r%
}}

\def\putAtriangle{\@ifnextchar <{\putAtrianglep}{\putAtrianglep
   <\arrowtypea`\arrowtypeb`\arrowtypec;\height>}}
\def\Atriangle{\@ifnextchar <{\Atrianglep}{\Atrianglep
   <\arrowtypea`\arrowtypeb`\arrowtypec;\height>}}

\def\Atrianglep<#1>[#2`#3`#4;#5`#6`#7]{{
\settriparms[#1]
\width=\height                         
\xext=\width                           
\yext=\height                          
\topadjust[#2``]
\botadjust[#3`#4`{#7}]
\multiply \xext by2 
\leftsladjust[#3`#2`{#5}]
\rightsladjust[#4`#2`{#6}]
\begin{picture}(\xext,\yext)(\xoff,\yoff)%
\putAtrianglep<\arrowtypea`\arrowtypeb`\arrowtypec;\height>%
(0,0)[#2`#3`#4;#5`#6`{#7}]%
\end{picture}%
}}

\def\putAtrianglepairp<#1>(#2)[#3;#4`#5`#6`#7`#8]{{
\settripairparms[#1]%
\setpos(#2)%
\settokens[#3]%
\puthmorphism(\xpos,\ypos)[\tokenb`\tokenc`{#7}]{\height}{\arrowtyped}b%
\advance\xpos by\height
\advance\ypos by\height
\putmorphism(\xpos,\ypos)(-1,-1)[\tokena``{#4}]{\height}{\arrowtypea}l%
\putvmorphism(\xpos,\ypos)[``{#5}]{\height}{\arrowtypeb}m%
\putmorphism(\xpos,\ypos)(1,-1)[``{#6}]{\height}{\arrowtypec}r%
}}

\def\putAtrianglepair{\@ifnextchar <{\putAtrianglepairp}{\putAtrianglepairp%
   <\arrowtypea`\arrowtypeb`\arrowtypec`\arrowtyped`\arrowtypee;\height>}}
\def\Atrianglepair{\@ifnextchar <{\Atrianglepairp}{\Atrianglepairp%
   <\arrowtypea`\arrowtypeb`\arrowtypec`\arrowtyped`\arrowtypee;\height>}}

\def\Atrianglepairp<#1>[#2;#3`#4`#5`#6`#7]{{%
\settripairparms[#1]%
\settokens[#2]%
\width=\height
\xext=\width
\yext=\height
\topadjust[\tokena``]%
\vertadjust[\tokenb`\tokenc`{#6}]
\tempcountd=\tempcounta                       
\vertadjust[\tokenc`\tokend`{#7}]
\ifnum\tempcounta<\tempcountd                 
\tempcounta=\tempcountd\fi                    
\advance \yext by\tempcounta                  
\advance \yoff by-\tempcounta                 %
\multiply \xext by2 
\leftsladjust[\tokenb`\tokena`{#3}]
\rightsladjust[\tokend`\tokena`{#5}]%
\begin{picture}(\xext,\yext)(\xoff,\yoff)%
\putAtrianglepairp
<\arrowtypea`\arrowtypeb`\arrowtypec`\arrowtyped`\arrowtypee;\height>%
(0,0)[#2;#3`#4`#5`#6`{#7}]%
\end{picture}%
}}

\def\putVtrianglep<#1>(#2,#3)[#4`#5`#6;#7`#8`#9]{{%
\settriparms[#1]%
\xpos=#2 \ypos=#3
\advance\ypos by\height
{\multiply\height by2
\puthmorphism(\xpos,\ypos)[#4`#5`{#7}]{\height}{\arrowtypea}a}%
\putmorphism(\xpos,\ypos)(1,-1)[`#6`{#8}]{\height}{\arrowtypeb}l%
\advance\xpos by\height
\advance\xpos by\height
\putmorphism(\xpos,\ypos)(-1,-1)[``{#9}]{\height}{\arrowtypec}r%
}}

\def\putVtriangle{\@ifnextchar <{\putVtrianglep}{\putVtrianglep
   <\arrowtypea`\arrowtypeb`\arrowtypec;\height>}}
\def\Vtriangle{\@ifnextchar <{\Vtrianglep}{\Vtrianglep
   <\arrowtypea`\arrowtypeb`\arrowtypec;\height>}}

\def\Vtrianglep<#1>[#2`#3`#4;#5`#6`#7]{{
\settriparms[#1]
\width=\height                         
\xext=\width                           
\yext=\height                          
\topadjust[#2`#3`{#5}]
\botadjust[#4``]
\multiply \xext by2 
\leftsladjust[#2`#3`{#6}]
\rightsladjust[#3`#4`{#7}]
\begin{picture}(\xext,\yext)(\xoff,\yoff)%
\putVtrianglep<\arrowtypea`\arrowtypeb`\arrowtypec;\height>%
(0,0)[#2`#3`#4;#5`#6`{#7}]%
\end{picture}%
}}

\def\putVtrianglepairp<#1>(#2)[#3;#4`#5`#6`#7`#8]{{
\settripairparms[#1]%
\setpos(#2)%
\settokens[#3]%
\advance\ypos by\height
\putmorphism(\xpos,\ypos)(1,-1)[`\tokend`{#6}]{\height}{\arrowtypec}l%
\puthmorphism(\xpos,\ypos)[\tokena`\tokenb`{#4}]{\height}{\arrowtypea}a%
\advance\xpos by\height
\putvmorphism(\xpos,\ypos)[``{#7}]{\height}{\arrowtyped}m%
\advance\xpos by\height
\putmorphism(\xpos,\ypos)(-1,-1)[``{#8}]{\height}{\arrowtypee}r%
}}

\def\putVtrianglepair{\@ifnextchar <{\putVtrianglepairp}{\putVtrianglepairp%
    <\arrowtypea`\arrowtypeb`\arrowtypec`\arrowtyped`\arrowtypee;\height>}}
\def\Vtrianglepair{\@ifnextchar <{\Vtrianglepairp}{\Vtrianglepairp%
    <\arrowtypea`\arrowtypeb`\arrowtypec`\arrowtyped`\arrowtypee;\height>}}

\def\Vtrianglepairp<#1>[#2;#3`#4`#5`#6`#7]{{%
\settripairparms[#1]%
\settokens[#2]
\xext=\height                  
\width=\height                 
\yext=\height                  
\vertadjust[\tokena`\tokenb`{#4}]
\tempcountd=\tempcounta        
\vertadjust[\tokenb`\tokenc`{#5}]
\ifnum\tempcounta<\tempcountd%
\tempcounta=\tempcountd\fi
\advance \yext by\tempcounta
\botadjust[\tokend``]%
\multiply \xext by2
\leftsladjust[\tokena`\tokend`{#6}]%
\rightsladjust[\tokenc`\tokend`{#7}]%
\begin{picture}(\xext,\yext)(\xoff,\yoff)%
\putVtrianglepairp
<\arrowtypea`\arrowtypeb`\arrowtypec`\arrowtyped`\arrowtypee;\height>%
(0,0)[#2;#3`#4`#5`#6`{#7}]%
\end{picture}%
}}

\def\putCtrianglep<#1>(#2,#3)[#4`#5`#6;#7`#8`#9]{{%
\settriparms[#1]%
\xpos=#2 \ypos=#3
\advance\ypos by\height
\putmorphism(\xpos,\ypos)(1,-1)[``{#9}]{\height}{\arrowtypec}l%
\advance\xpos by\height
\advance\ypos by\height
\putmorphism(\xpos,\ypos)(-1,-1)[#4`#5`{#7}]{\height}{\arrowtypea}l%
{\multiply\height by 2
\putvmorphism(\xpos,\ypos)[`#6`{#8}]{\height}{\arrowtypeb}r}%
}}

\def\putCtriangle{\@ifnextchar <{\putCtrianglep}{\putCtrianglep
    <\arrowtypea`\arrowtypeb`\arrowtypec;\height>}}
\def\Ctriangle{\@ifnextchar <{\Ctrianglep}{\Ctrianglep
    <\arrowtypea`\arrowtypeb`\arrowtypec;\height>}}

\def\Ctrianglep<#1>[#2`#3`#4;#5`#6`#7]{{
\settriparms[#1]
\width=\height                          
\xext=\width                            
\yext=\height                           
\multiply \yext by2 
\topadjust[#2``]
\botadjust[#4``]
\sladjust[#3`#2`{#5}]{\width}
\tempcountd=\tempcounta                 
\sladjust[#3`#4`{#7}]{\width}
\ifnum \tempcounta<\tempcountd          
\tempcounta=\tempcountd\fi              
\advance \xext by\tempcounta            
\advance \xoff by-\tempcounta           %
\rightadjust[#2`#4`{#6}]
\begin{picture}(\xext,\yext)(\xoff,\yoff)%
\putCtrianglep<\arrowtypea`\arrowtypeb`\arrowtypec;\height>%
(0,0)[#2`#3`#4;#5`#6`{#7}]%
\end{picture}%
}}

\def\putDtrianglep<#1>(#2,#3)[#4`#5`#6;#7`#8`#9]{{%
\settriparms[#1]%
\xpos=#2 \ypos=#3
\advance\xpos by\height \advance\ypos by\height
\putmorphism(\xpos,\ypos)(-1,-1)[``{#9}]{\height}{\arrowtypec}r%
\advance\xpos by-\height \advance\ypos by\height
\putmorphism(\xpos,\ypos)(1,-1)[`#5`{#8}]{\height}{\arrowtypeb}r%
{\multiply\height by 2
\putvmorphism(\xpos,\ypos)[#4`#6`{#7}]{\height}{\arrowtypea}l}%
}}

\def\putDtriangle{\@ifnextchar <{\putDtrianglep}{\putDtrianglep
    <\arrowtypea`\arrowtypeb`\arrowtypec;\height>}}
\def\Dtriangle{\@ifnextchar <{\Dtrianglep}{\Dtrianglep
   <\arrowtypea`\arrowtypeb`\arrowtypec;\height>}}

\def\Dtrianglep<#1>[#2`#3`#4;#5`#6`#7]{{
\settriparms[#1]
\width=\height                         
\xext=\height                          
\yext=\height                          
\multiply \yext by2 
\topadjust[#2``]
\botadjust[#4``]
\leftadjust[#2`#4`{#5}]
\sladjust[#3`#2`{#5}]{\height}
\tempcountd=\tempcountd                
\sladjust[#3`#4`{#7}]{\height}
\ifnum \tempcounta<\tempcountd         
\tempcounta=\tempcountd\fi             
\advance \xext by\tempcounta           %
\begin{picture}(\xext,\yext)(\xoff,\yoff)
\putDtrianglep<\arrowtypea`\arrowtypeb`\arrowtypec;\height>%
(0,0)[#2`#3`#4;#5`#6`{#7}]%
\end{picture}%
}}

\def\setrecparms[#1`#2]{\width=#1 \height=#2}%
\def\recursep<#1`#2>[#3;#4`#5`#6`#7`#8]{{%
\width=#1 \height=#2
\settokens[#3]
\settowidth{\tempdimen}{$\tokena$}
\ifdim\tempdimen=0pt
  \savebox{\tempboxa}{\hbox{$\tokenb$}}%
  \savebox{\tempboxb}{\hbox{$\tokend$}}%
  \savebox{\tempboxc}{\hbox{$#6$}}%
\else
  \savebox{\tempboxa}{\hbox{$\hbox{$\tokena$}\times\hbox{$\tokenb$}$}}%
  \savebox{\tempboxb}{\hbox{$\hbox{$\tokena$}\times\hbox{$\tokend$}$}}%
  \savebox{\tempboxc}{\hbox{$\hbox{$\tokena$}\times\hbox{$#6$}$}}%
\fi
\ypos=\height
\divide\ypos by 2
\xpos=\ypos
\advance\xpos by \width
\xext=\xpos \yext=\height
\topadjust[#3`\usebox{\tempboxa}`{#4}]%
\botadjust[#5`\usebox{\tempboxb}`{#8}]%
\sladjust[\tokenc`\tokenb`{#5}]{\ypos}%
\tempcountd=\tempcounta
\sladjust[\tokenc`\tokend`{#5}]{\ypos}%
\ifnum \tempcounta<\tempcountd
\tempcounta=\tempcountd\fi
\advance \xext by\tempcounta
\advance \xoff by-\tempcounta
\rightadjust[\usebox{\tempboxa}`\usebox{\tempboxb}`\usebox{\tempboxc}]%
\bfig
\putCtrianglep<-1`1`1;\ypos>(0,0)[`\tokenc`;#5`#6`{#7}]%
\puthmorphism(\ypos,0)[\tokend`\usebox{\tempboxb}`{#8}]{\width}{-1}b%
\puthmorphism(\ypos,\height)[\tokenb`\usebox{\tempboxa}`{#4}]{\width}{-1}a%
\advance\ypos by \width
\putvmorphism(\ypos,\height)[``\usebox{\tempboxc}]{\height}1r%
\efig
}}

\def\recurse{\@ifnextchar <{\recursep}{\recursep<\width`\height>}}

\def\puttwohmorphisms(#1,#2)[#3`#4;#5`#6]#7#8#9{{%
\puthmorphism(#1,#2)[#3`#4`]{#7}0a
\ypos=#2
\advance\ypos by 20
\puthmorphism(#1,\ypos)[\phantom{#3}`\phantom{#4}`#5]{#7}{#8}a
\advance\ypos by -40
\puthmorphism(#1,\ypos)[\phantom{#3}`\phantom{#4}`#6]{#7}{#9}b
}}

\def\puttwovmorphisms(#1,#2)[#3`#4;#5`#6]#7#8#9{{%
\putvmorphism(#1,#2)[#3`#4`]{#7}0a
\xpos=#1
\advance\xpos by -20
\putvmorphism(\xpos,#2)[\phantom{#3}`\phantom{#4}`#5]{#7}{#8}l
\advance\xpos by 40
\putvmorphism(\xpos,#2)[\phantom{#3}`\phantom{#4}`#6]{#7}{#9}r
}}

\def\puthcoequalizer(#1)[#2`#3`#4;#5`#6`#7]#8#9{{%
\setpos(#1)%
\puttwohmorphisms(\xpos,\ypos)[#2`#3;#5`#6]{#8}11%
\advance\xpos by #8
\puthmorphism(\xpos,\ypos)[\phantom{#3}`#4`#7]{#8}1{#9}
}}

\def\putvcoequalizer(#1)[#2`#3`#4;#5`#6`#7]#8#9{{%
\setpos(#1)%
\puttwovmorphisms(\xpos,\ypos)[#2`#3;#5`#6]{#8}11%
\advance\ypos by -#8
\putvmorphism(\xpos,\ypos)[\phantom{#3}`#4`#7]{#8}1{#9}
}}

\def\putthreehmorphisms(#1)[#2`#3;#4`#5`#6]#7(#8)#9{{%
\setpos(#1) \settypes(#8)
\if a#9 %
     \vertsize{\tempcounta}{#5}%
     \vertsize{\tempcountb}{#6}%
     \ifnum \tempcounta<\tempcountb \tempcounta=\tempcountb \fi
\else
     \vertsize{\tempcounta}{#4}%
     \vertsize{\tempcountb}{#5}%
     \ifnum \tempcounta<\tempcountb \tempcounta=\tempcountb \fi
\fi
\advance \tempcounta by 60
\puthmorphism(\xpos,\ypos)[#2`#3`#5]{#7}{\arrowtypeb}{#9}
\advance\ypos by \tempcounta
\puthmorphism(\xpos,\ypos)[\phantom{#2}`\phantom{#3}`#4]{#7}{\arrowtypea}{#9}
\advance\ypos by -\tempcounta \advance\ypos by -\tempcounta
\puthmorphism(\xpos,\ypos)[\phantom{#2}`\phantom{#3}`#6]{#7}{\arrowtypec}{#9}
}}

\def\putarc(#1,#2)[#3`#4`#5]#6#7#8{{%
\xpos #1
\ypos #2
\width #6
\arrowlength #6
\putbox(\xpos,\ypos){#3\vphantom{#4}}%
{\advance \xpos by\arrowlength
\putbox(\xpos,\ypos){\vphantom{#3}#4}}%
\horsize{\tempcounta}{#3}%
\horsize{\tempcountb}{#4}%
\divide \tempcounta by2
\divide \tempcountb by2
\advance \tempcounta by30
\advance \tempcountb by30
\advance \xpos by\tempcounta
\advance \arrowlength by-\tempcounta
\advance \arrowlength by-\tempcountb
\halflength=\arrowlength \divide\halflength by 2
\divide\arrowlength by 5
\put(\xpos,\ypos){\bezier{\arrowlength}(0,0)(50,50)(\halflength,50)}
\ifnum #7=-1 \put(\xpos,\ypos){\vector(-3,-2)0} \fi
\advance\xpos by \halflength
\put(\xpos,\ypos){\xpos=\halflength \advance\xpos by -50
   \bezier{\arrowlength}(0,50)(\xpos,50)(\halflength,0)}
\ifnum #7=1 {\advance \xpos by
   \halflength \put(\xpos,\ypos){\vector(3,-2)0}} \fi
\advance\ypos by 50
\vertsize{\tempcounta}{#5}%
\divide\tempcounta by2
\advance \tempcounta by20
\if a#8 %
   \advance \ypos by\tempcounta
   \putbox(\xpos,\ypos){#5}%
\else
   \advance \ypos by-\tempcounta
   \putbox(\xpos,\ypos){#5}%
\fi
}}

\makeatother

\usepackage{amsthm}
\usepackage{dsfont}
\usepackage{stmaryrd}
\newtheorem{theorem}{Theorem}[section]
\newtheorem{lemma}[theorem]{Lemma}

\makeindex \makeglossary
\begin{document}

\sloppy

\newcommand{\nl}{\hspace{2cm}\\ }

\def\nec{\Box}
\def\pos{\Diamond}
\def\diam{{\tiny\Diamond}}

\def\lc{\lceil}
\def\rc{\rceil}
\def\lf{\lfloor}
\def\rf{\rfloor}
\def\lk{\langle}
\def\rk{\rangle}
\def\blk{\dot{\langle\!\!\langle}}
\def\brk{\dot{\rangle\!\!\rangle}}

\newcommand{\pa}{\parallel}
\newcommand{\lra}{\longrightarrow}
\newcommand{\hra}{\hookrightarrow}
\newcommand{\hla}{\hookleftarrow}
\newcommand{\ra}{\rightarrow}
\newcommand{\la}{\leftarrow}
\newcommand{\lla}{\longleftarrow}
\newcommand{\da}{\downarrow}
\newcommand{\ua}{\uparrow}
\newcommand{\dA}{\downarrow\!\!\!^\bullet}
\newcommand{\uA}{\uparrow\!\!\!_\bullet}
\newcommand{\Da}{\Downarrow}
\newcommand{\DA}{\Downarrow\!\!\!^\bullet}
\newcommand{\UA}{\Uparrow\!\!\!_\bullet}
\newcommand{\Ua}{\Uparrow}
\newcommand{\Lra}{\Longrightarrow}
\newcommand{\Ra}{\Rightarrow}
\newcommand{\Lla}{\Longleftarrow}
\newcommand{\La}{\Leftarrow}
\newcommand{\nperp}{\perp\!\!\!\!\!\setminus\;\;}
\newcommand{\pq}{\preceq}

\newcommand{\lms}{\longmapsto}
\newcommand{\ms}{\mapsto}
\newcommand{\subseteqnot}{\subseteq\hskip-4 mm_\not\hskip3 mm}

\def\o{{\omega}}

\def\bA{{\bf A}}
\def\bEM{{EM}}
\def\bM{{\bf M}}
\def\bN{{\bf N}}
\def\bC{{\bf C}}
\def\bI{{\bf I}}
\def\bK{{K}}
\def\bL{{\bf L}}
\def\bT{{\bf T}}
\def\bS{{\bf S}}
\def\bD{{\bf D}}
\def\bB{{\bf B}}
\def\bW{{\bf W}}
\def\bP{{\bf P}}
\def\bX{{\bf X}}
\def\bY{{\bf Y}}
\def\ba{{\bf a}}
\def\bb{{\bf b}}
\def\bc{{\bf c}}
\def\bd{{\bf d}}
\def\bh{{\bf h}}
\def\bi{{\bf i}}
\def\bj{{\bf j}}
\def\bk{{\bf k}}
\def\bm{{\bf m}}
\def\bn{{\bf n}}
\def\bp{{\bf p}}
\def\bq{{\bf q}}
\def\be{{\bf e}}
\def\br{{\bf r}}
\def\bi{{\bf i}}
\def\bs{{\bf s}}
\def\bt{{\bf t}}
\def\jeden{{\bf 1}}
\def\dwa{{\bf 2}}
\def\trzy{{\bf 3}}

\def\cB{{\cal B}}
\def\cA{{\cal A}}
\def\cC{{\cal C}}
\def\cD{{\cal D}}
\def\cE{{\cal E}}
\def\cEM{{\cal EM}}
\def\cF{{\cal F}}
\def\cG{{\cal G}}
\def\cI{{\cal I}}
\def\cJ{{\cal J}}
\def\cK{{\cal K}}
\def\cL{{\cal L}}
\def\cN{{\cal N}}
\def\cM{{\cal M}}
\def\cO{{\cal O}}
\def\cP{{\cal P}}
\def\cQ{{\cal Q}}
\def\cR{{\cal R}}
\def\cS{{\cal S}}
\def\cT{{\cal T}}
\def\cU{{\cal U}}
\def\cV{{\cal V}}
\def\cW{{\cal W}}
\def\cX{{\cal X}}
\def\cY{{\cal Y}}


\def\Mnd{{\bf Mnd}}
\def\AMnd{{\bf AnMnd}}
\def\An{{\bf An}}
\def\San{{\bf San}}
\def\PMnd{{\bf PolyMnd}}
\def\SanMnd{{\bf SanMnd}}
\def\RiMnd{{\bf RiMnd}}
\def\End{{\bf End}}

\def\ET{\bf ET}
\def\RegET{\bf RegET}
\def\RET{\bf RegET}
\def\LrET{\bf LrET}
\def\RiET{\bf RiET}
\def\SregET{\bf SregET}

\def\LT{\bf LT}
\def\RegLT{\bf RegLT}
\def\ALT{\bf AnLT}
\def\RiLT{\bf RiLT}

\def\FOp{\bf FOp}
\def\RegOp{\bf RegOp}
\def\SOp{\bf SOp}
\def\RiOp{\bf RiOp}

\def\bCat{{{\bf Cat}}}
\def\MonCat{{{\bf MonCat}}}
\def\Mon{{{\bf Mon}}}
\def\Cat{{{\bf Cat}}}

\def\F{\mathds{F}}
\def\S{\mathds{S}}
\def\I{\mathds{I}}
\def\B{\mathds{B}}

\def\V{\mathds{V}}
\def\W{\mathds{W}}
\def\M{\mathds{M}}
\def\N{\mathds{N}}

\def\Op{{\cal O}p}

\def\Vb{\bar{\mathds{V}}}
\def\Wb{\bar{\mathds{W}}}

\def\Sym{{\cal S}ym}

\def\P{{\cal P}}
\def\Q{{\cal Q}}

\pagenumbering{arabic} \setcounter{page}{1}

\title{\bf\Large \Large Rigidity is undecidable}
\author{Miko\l aj Bojanczyk, Stanis\l aw Szawiel,\\ Marek Zawadowski}

\maketitle
\begin{abstract}  We show that the problem `whether a finite set of regular-linear axioms defines a rigid theory' is undecidable.
\end{abstract}
{\em 2010 Mathematical Subject Classification} 03D35, 03C05, 03G30, 18C10, 18C15

{\em Keywords:} Equational theory, interpretation, undecidable problem, word problem for monoid
\section{Introduction}

In \cite{SZ} it was shown that the category of polynomial monads is equivalent to the category of rigid equational theories, solving a problem stated in \cite{CJ2}. A linear-regular theory is an equational theory that has as a set of axioms equations of terms $s=t$ such that the variables occurring in $s$ and $t$ are the same and each of them occurs  once. For example the theory of monoids, commutative monoids, and monoids with anti-involution are linear-regular. Recall that the theory of monoids with anti-involution has three function symbols e, i, m of arity 0, 1, 2, respectively, contains the usual axioms for monoids, and additionally the equations $i(i(x_1)=x_1$ and $i(m(x_1,x_2))=m(i(x_2),i(x_1))$. A linear-regular theory is rigid if and only if for any term $t(x_1,\ldots,x_n)$ and permutation $\sigma\in S_n$, if $T\vdash t(x_1,\ldots,x_n)=t(x_{\sigma(1)},\ldots,x_{\sigma(n)})$ then $\sigma$ is an identity permutation. In other words, $T$ is rigid if any (proper) permutation of variables changes the meanings of terms in $T$. For example, the theories of monoids and monoids with anti-involution are rigid but the theory of commutative monoids is not as it contains the axiom $m(x_1,x_2)=m(x_2,x_1)$. Rigidity refers to provability in $T$ and hence it is a global property concerning a linear-regular theories. In this paper we show that the problem whether an equational theory $T$ with finite set of linear-regular axioms is rigid, is undecidable.

\section{Preliminaries}

When dealing with equational theories we follow mostly the terminology of \cite{BN}. However, we want to specify what variables might occur in a term, and for this reason we deal with terms in context. We call an `equation' what in \cite{BN} is called an `identity'.

By an {\em equational theory} we mean a pair of sets $T=(L,A)$, $L=\bigcup_{n\in \o} L_n$ and $L_n$ is the set of $n$-ary operations.
The sets of operations of different arities are disjoint.
The set $\cT r(L,\vec{x}^n)$ of terms of $L$ in context $\vec{x}^n=\lk x_1, \ldots, x_n \rk $
is the usual set of terms over $L$ built with the help of variables from $\vec{x}^n$.
We write $t:\vec{x}^n$ for the term $t$ in context $\vec{x}^n$. Thus all the variables occurring in $t$ are among those in $\vec{x}^n$.
The set $A$ is a set of equations in context $t=s : \vec{x}^n$, i.e.  both $t:\vec{x}^n$ and $s:\vec{x}^n$ are terms in context, $n\in\o$.

If we do not specify explicitly the context of a term then we mean that the context consists of variables explicitly occurring in the term. As in \cite{BN} we often think of a term as a tree labeled by functions symbols and variables. A derivation consists of a finite number of rewrite steps. One rewrite step replaces a part of a term tree that matches a substitution of one side of an equation in $T$ by the same substitution of the other side of that equation. For details see \cite{BN} definition  3.1.8. When possible, a simple derivation will be presented as a sequence of equations.

A morphism of equational theories, an {\em interpretation},  $I : T \ra T'=(L',A')$, is a set of functions $I_n : L_n \ra \cT r(L',\vec{x}^n)$ for $n\in \o$. Moreover, we require that $I$ preserves the equations, i.e. for any $t=s : \vec{x}^n$ in $A$ we have
\[ A' \vdash \bar{I}(t)=\bar{I}(s) : \vec{x}^n\]
where $A'\vdash$ (or $T'\vdash$) is the provability in the equational logic from axioms in the set of axioms $A'$ (or theory $T'$).  $\bar{I}$ is the extension of
$I_n$'s to functions  $\bar{I}_n : \cT r(L,\vec{x}^n) \ra \cT r(L',\vec{x}^n)$ for $n\in \o$ as follows. We usually drop index $n$ in $\bar{I}_n$.

\[ \bar{I}(x_i:\vec{x}^n) = x_i:\vec{x}^n \] for $i=1, \ldots, n$, $n\in\o$ and

\[ \bar{I}(f(t_1,\ldots, t_k) :\vec{x}^n) = I(f)(x_1\setminus \bar{I}(t_1),\ldots, x_k\setminus \bar{I}(t_k)):\vec{x}^n \]
for $f\in L_k$ and $t_i\in \cT r(L,\vec{x}^n)$ for $i=1,\ldots, k$. On the right-hand side, we have a simultaneous substitution of terms $t_i$'s for variables $x_i$'s.  We identify two such interpretations $I$ and  $I' : (L,A) \ra (L',A')$
iff they interpret all function symbols as provably equivalent terms, i.e.
\[ A' \vdash I(f)=I'(f) : \vec{x}^n\]
for any $n\in\o$ and $f\in L_n$. An interpretation $I : T \ra T'$ is {\em conservative} iff for any equation in context $s=t:\vec{x}^n$ in $T$
if  $T' \vdash \bar{I}(s)=\bar{I}(t) : \vec{x}^n$ then  $T \vdash s=t : \vec{x}^n$.

A term in context $t:\vec{x}^n$ is {\em linear-regular} if every variable in $\vec{x}^n$ occurs in $t$ exactly once.
An equation $s=t:\vec{x}^n$ is {\em linear-regular} iff both $s:\vec{x}^n$ and
$t:\vec{x}^n$ are linear-regular terms in contexts.

A {\em simple $\phi$-substitution of a term in context}  $t:\vec{x}^n$ along a function $\phi:(n]\ra (k]$ is a term in context denoted $\phi\cdot t : \vec{x}^k$ such that every occurrence of the variable $x_i$ is replaced by the occurrence of $x_{\phi(i)}$.

An equational theory $T=(L,A)$ is a {\em linear-regular theory} iff all the consequences of $T$  are consequences of linear-regular consequences of $T$.
An interpretation is a {\em linear-regular interpretation} iff it interprets function symbols as linear-regular terms.

A theory $T=(L,A)$ is a {\em rigid theory} iff it is linear-regular and for any linear-regular term in context $t:\vec{x}^n$ whenever  $A\vdash t=\sigma \cdot t:\vec{x}^n$ then $\sigma$ is the identity permutation. $\tau \cdot t$ is the simple $\sigma$-substitution of a term in context  $t:\vec{x}^n$ along a permutation $\sigma\in S_n$.

The definitions of both linear-regular and rigid theories are such to make sure that if a theory is isomorphic to a linear-regular (rigid) theory then it is also linear-regular (rigid).

\section{Main result}

If we find a linear-regular set of axioms of an equational theory $T$ we can be sure that $T$ is linear-regular, (cf. \cite{SZ}). However, it is not so easy to decide whether a given theory is rigid. The main result of this paper says that even if we restrict ourselves to finitely axiomatizable linear-regular theories it is still undecidable whether such theories are rigid or not.

A term $t(x_1,\ldots, x_n)$ is {\em flabby} in $T$ if it is linear-regular in variables $x_1,\ldots, x_n$ such that
\begin{equation} \label{flabby}
T\vdash t(x_1,\ldots, x_n)=t(x_{\sigma(1)},\ldots, x_{\sigma(n)})
\end{equation}
for a non-identity permutation  $\sigma\in S_n$. A theory is rigid iff it does not contain flabby terms.

\begin{theorem}\label{Undecidable}
The problem whether an equational theory $T=(L,A)$ in finite language $L$
with a finite set of linear-regular axioms $A$ is rigid is undecidable.
\end{theorem}

{\em Proof.} The word problem for monoids is undecidable; (cf. \cite{M}, \cite{P}). We shall show that it reduces to our problem. Below we sketch the construction of the reduction and an argument showing that it is indeed a reduction. Then in a series of Lemmas proved in the remaining part of the paper we shall make the sketched construction and argument more precise.

First we define a simple theory $T_0$ that is rigid, (cf. Lemma \ref{T0rigid}). For an arbitrary instance of the word problem for monoids,
\begin{equation} \label{word-problem} \bigwedge_{i\in n} u_i=v_i \vdash u=v
\end{equation}
where $u_i,v_i,u,v$ are words over a finite alphabet, we will define a theory $T$ such that $T$ is rigid iff (\ref{word-problem}) does not hold.

An easy argument shows that if (\ref{word-problem}) holds then there is an obvious flabby term in $T$ and hence $T$ is not rigid, (cf. Lemma \ref{not rigid}).

Next we define a linear-regular interpretation $I:T_0\ra T$ which is conservative iff (\ref{word-problem}) does not hold, (cf. Lemma \ref{cons interpretation}).
The terms in the image of $\bar{I}: \cT r(T_0)\ra  \cT r(T)$ are called special and the set of special terms is denoted by $\cS p(T)$. We construct a function
\[ \widehat{(-)} : \cT r(T)\lra \cS p(T) \]
sending all terms of $T$ to the special terms such that
\begin{enumerate}
  \item $\widehat{(-)}$ is onto;
  \item $\widehat{\bar{I}(s)}=\bar{I}(s)$, for any  $s\in \cT r(T_0)$;
  \item for $t\in \cT r(T)$, the variables occurring in both terms $t$ and $\hat{t}$ are the same and they occur in the same order;
  \item for $t,t'\in \cT r(T)$, if $T\vdash t=t'$ then $T\vdash \hat{t}=\hat{t'}$;
\end{enumerate}
(cf. Lemma \ref{hat_function}).

Having established the above, to get a contradiction, we shall assume that (\ref{word-problem}) does not hold but $T$ is still not rigid.
Let $t(x_1,\ldots , x_n)$ be a flabby term in $T$ and $\sigma\in S_n$ such that (\ref{flabby}) holds. Then, by Lemma \ref{hat_function},
\begin{equation} \label{flabby1}
T\vdash \hat{t}(x_1,\ldots, x_n)=\hat{t(}x_{\sigma(1)},\ldots, x_{\sigma(n)})
\end{equation}
holds.  As $\bar{I}$ is onto there is a term $s(x_1,\ldots, x_n)$  in $T_0$ such that $\bar{I}(s)(x_1,\ldots, x_n)=  \hat{t}(x_1,\ldots, x_n)$. Thus
$$ T\vdash \bar{I}(s)(x_1,\ldots, x_n)=\bar{I}(s)(x_{\sigma(1)},\ldots, x_{\sigma(n)})$$
and since $I$ is conservative
$$ T_0\vdash s(x_1,\ldots, x_n)=s(x_{\sigma(1)},\ldots, x_{\sigma(n)}).$$
But this mean that $s$ is a flabby term in $T_0$, contradicting rigidity of $T_0$. This ends the proof of the theorem.
$\boxempty$

Now, we fix for the rest of the paper the theory $T$ constructed as in the (sketch of) proof of Theorem \ref{Undecidable} and we fill the details of the above argument.

The theory $T_0$ contains three binary symbols $l,r,m$ and one equation
\begin{equation} \label{identity3}
l(x_1,x_2)=r(x_2,x_1)
\end{equation}
We have

\begin{lemma}\label{T0rigid}
$T_0$ is a rigid theory.
\end{lemma}
{\em Proof.}

The theory $T_0$ is equivalent (in fact isomorphic) to a theory that has two binary function symbols and no equations. Thus it contains no non-trivial equations. In particular it is rigid.
$\boxempty$

Let us fix an instance of the word problem for monoids. Let  $u_i,v_i,u,v$ words over the alphabet $G=\{ g_1,\ldots, g_n \}$, for $i\in m$. The problem is to decide whether (\ref{word-problem}) holds true. We define an equational theory $T$ corresponding to this problem. The alphabet of $T$ consists of unary symbols from $G$ and additionally one unary symbol $\alpha$ and one binary symbol $m$. If $w=g_{k_1}\ldots g_{k_m}$ is a word over $G$ then $w(x)$ denotes the corresponding term $g_{k_1}\ldots g_{k_m}(x)$ of $T$. The axioms of $T$ are
\begin{equation} \label{identity1}
u_i(x_1)=v_i(x_1) \;\;\; {\rm for}\;\;\; i\in m
\end{equation}
and moreover
\begin{equation} \label{identity2}
 m(u\alpha (x_1),x_2)= m(v\alpha (x_2),x_1)
\end{equation}

The following Lemma makes simple but useful observations concerning the derivations in theory $T$.

\begin{lemma}\label{derivarions_in_T}
\begin{enumerate}
\item For any two words $w_1$, $w_2$ over $G$ we have
\begin{equation} \label{consequence}
 \bigwedge_{i\in m} u_i=v_i \vdash w_1=w_2\;\;\; {\rm iff}\;\;\; T \vdash w_1(x_1)=w_2(x_1)
\end{equation}
where $\vdash$ on the left is the consequence relation in the theory of monoids.
  \item The symbol $\alpha$ does not take part in any rewrite step over $T$ concerning unary symbols.
  \item Each rewrite step over $T$ concerns only unary symbols or it is performed on a subterm with the root labeled $m$. In particular, no derivation changes the number of symbols $m$. $\boxempty$
\end{enumerate}
\end{lemma}

{\em Remark.} Last property says that in the derivations in $T$ we can trace the identity of each symbol $m$. We are going to use it when arguing about derivations.

\begin{lemma}\label{not rigid}
If (\ref{word-problem}) holds then $T$ is not rigid.
\end{lemma}
{\em Proof.} Let $t(x_1,x_2)=m(u\alpha (x_1),x_2)$. Then, using (\ref{identity2}), (\ref{word-problem}), and \ref{consequence}, we have in $T$
\[ t(x_1,x_2) = m(u\alpha (x_1),x_2)= m(v\alpha (x_2),x_1)= m(u\alpha (x_2),x_1)=t(x_2,x_1) \]
i.e. $t$ is flabby in $T$, and $T$ is not rigid.
$\boxempty$
\vskip 2mm

Now we define a linear-regular interpretation $I:T_0\ra T$ as follows
\[ I(l)=m(u\alpha (x_1),x_2),\;\;  I(r)=m(v\alpha (x_1),x_2),\;\; I(m)=m(x_1,x_2). \]

\begin{lemma}\label{cons interpretation}
$I:T_0\ra T$ is a linear-regular interpretation. It is conservative iff (\ref{word-problem}) does not hold.
\end{lemma}
{\em Proof.} We have in $T$
\[ \bar{I}(l)(x_1,x_2)=m(u\alpha (x_1),x_2)= m(v\alpha (x_2),x_1)=\bar{I}(l)(x_2,x_1) \]
and hence $I$ is an interpretation.

If (\ref{word-problem}) holds then we have  in $T$
\[ \bar{I}(l(x_1,x_2))= m(u\alpha (x_1),x_2)= m(v\alpha (x_2),x_1)=  \bar{I}(r(x_1,x_2)) \]
But clearly $T_0\not\vdash l(x_1,x_2)=r(x_1,x_2)$. So $I$ is not conservative.

Now, we assume that $T\not\vdash u(x)=v(x)$ and we shall show that $I$ is conservative.  Let $s$, $s'$ be two terms in $T_0$ such that
\[  T\vdash\bar{I}(s)=\bar{I}(s'). \]
First, we want to show that the above equality can be deduced without use of the equations (\ref{identity1}). Let $D$ be a derivation of $\bar{I}(s)=\bar{I}(s')$ in $T$ that contains minimal number of applications of equations (\ref{identity1}). If $D$ does not use (\ref{identity1}), we are done. If the equation (\ref{identity1}) is used in $D$, it is used to either part of the string of unary symbols $u$ or $v$ of a subterm $m(u\alpha(t_1),t_2)$ or $m(v\alpha(t_1),t_2)$, respectively. Suppose the first rewrite step using the equation (\ref{identity1}) in the derivation $D$ is applied to the subterm $m(u\alpha(t_1),t_2)$ rewriting it to some other subterm $m(u'\alpha(t_1),t_2)$. The rewrite steps concerning the subterm with `this occurrence' of $m$ as the root symbol will concern the subterms $u'$ and $t_1$, $t_2$ parts only and possibly $m$ but only if $u'$ will be rewritten to either $u$ or $v$. By assumption, $u$ cannot be rewritten to $v$, so it can only be rewritten back to $u$. In fact, as at the end of the derivation we get a term of form $\bar{I}(s')$ (with all strings of unary symbols from $G$ being equal either $u$ or $v$), $u'$ has to be eventually rewritten back to $u$. But this means that we can shorten the derivation $D$ by eliminating all those rewrite steps from $u$ to $u'$ and back to $u$ again. As this contradicts the minimality of $D$, we can assume that $D$ contains only rewrite steps that use the equation (\ref{identity2}). But then the derivation $D$ of $\bar{I}(s)=\bar{I}(s')$ in $T$ can be used to build a derivation $D'$ of  $s=s'$ in $T_0$. We need to change the rewrite steps using the equation (\ref{identity2}) in $D$ to rewrite steps in the corresponding  positions using the equality (\ref{identity3}) in $D'$. Thus  $T_0\vdash s=s'$. Since terms $s$, $s'$ where arbitrary,  $I$ is conservative. $\boxempty$
\vskip 2mm

Special terms of $T$ are terms in the image of the function $\bar{I}: \cT r (T_0) \ra \cT r(T)$. The set of special terms is denoted by $\cS p(T)$.
The function
\[ \widehat{(-)} : \cT r (T) \lra \cS p(T) \]
is defined, for $t=t(x_1,\ldots, x_n)\in\cT r (T)$  as follows
\[ \widehat{t}\;\;= \;\; \left\{ \begin{array}{ll}
                x_i   & \mbox{ if } t=x_i \\
                \widehat{t'}   & \mbox{ if } t=g(t')\mbox{ where } g\in G\cup \{ \alpha\} \\
                m(u\alpha (\widehat{t_1}),\widehat{t_2}) & \mbox{ if } t=m(w\alpha (t_1),t_2)\mbox{ and } T\vdash u(x)=w(x) \\
                m(v\alpha (\widehat{t_1}),\widehat{t_2}) & \mbox{ if } t=m(w\alpha (t_1),t_2) \mbox{ and } T\vdash v(x)=w(x)\mbox{ and not } (\ref{word-problem}) \\
                m(\widehat{t_1},\widehat{t_2}) & \mbox{ if } t=m(t_1,t_2), \mbox{ and none of the above applies}.
                                    \end{array}
                            \right. \]
The following Lemma lists some properties of $\widehat{(-)}$ that were used in the proof of the main theorem.

\begin{lemma}\label{hat_function} We have
\begin{enumerate}
  \item $\widehat{(-)}$ is onto;
  \item $\widehat{\bar{I}(s)}=\bar{I}(s)$, for any term $s$ of $T_0$;
  \item for any term $t$ of $T$, the variables in terms $t$ and $\widehat{t}$ are the same and they occur in the same order; $t : \vec{x}^n$ is a linear-regular term iff $\widehat{t}: \vec{x}^n$ is;
  \item if $T \vdash t=t'$ then  $T \vdash \widehat{t}=\widehat{t'}$,  for any terms $t$, $t'$ in $T$.
\end{enumerate}
\end{lemma}
{\em Proof.} 1. and 2. is obvious.

To show 3. one can verify by induction on the construction of terms that no clause in the definition of $\widehat{(-)}$ changes the variables or their order.

We shall show 4. by induction on the complexity of the term $t$. If $t$ is a variable then the thesis is obvious.

If $t=w(m(t_1,t_2))$  where $w$ is a (non-empty) sequence of unary symbols of $T$ then, as no derivation changes the number of symbols in terms, $t'=w'(m(t'_1,t'_2))$  where $w'$ is a sequence of unary symbols of $T$. The derivation $D$ from $t$ to $t'$ consists of steps that either change unary symbols over the first $m$ in the term using equations (\ref{identity1}) or does not involve those symbols at all. Thus if we drop from the derivation $D$ all the rewrite steps that change symbols over the first $m$ the resulting derivation proves $w(m(t_1,t_2))=w(m(t'_1,t'_2))$ never using symbols from $w$. Thus the same derivation proves also $m(t_1,t_2)=m(t'_1,t'_2)$. Using inductive hypothesis we get
\[ \widehat{t} = \widehat{m(t_1,t_2)} =\widehat{m(t'_1,t'_2)} = \widehat{t'} \]

If $t=m(t_1,t_2)$ and  $t'=z(m(t'_1,t'_2))$, then by the above we can assume that the sequence of the unary symbols $z$ is empty. We have to consider three cases concerning the form of the term $t_1$:
\begin{enumerate}
  \item $t_1= w(\alpha (s))$ and $T\vdash u(x)=w(x)$;
  \item $t_1= w(\alpha (s))$ and $T\vdash v(x)=w(x)$;
  \item $t_1$ is not in the above form.
\end{enumerate}
As the cases 1. and 2. are similar we shall consider cases 1. and 3 only.

We start with Case 3, as it is much simpler.
In that case to the leading symbol $m$ in term $t$ the rule (\ref{identity2}) is never applied. Thus all the derivations of $T\vdash t=t'$ can be split into two separate derivations, one for $T\vdash t_1=t'_1$ and one for $T\vdash t_2=t'_2$. Thus again using inductive hypothesis we get
\[ \widehat{t} = m(\widehat{t_1},\widehat{t_2}) =m(\widehat{t'_1},\widehat{t'_2}) = \widehat{t'} \;\;\boxempty\]

It remains to consider Case 1. The term $t$ looks as follows
  \begin{center} \xext=500 \yext=600
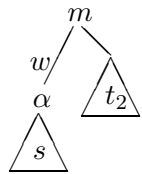
\begin{picture}(\xext,\yext)(\xoff,\yoff)
 \put(0,0){\line(1,0){200}}
  \put(0,0){\line(1,2){100}}
  \put(200,0){\line(-1,2){100}}
   \put(80,40){$s$}
 \put(80,220){$\alpha$}
 \put(120,300){\line(1,2){100}}
 \put(70,350){$w$}
   \put(200,520){$m$}
 \put(250,500){\line(1,-1){100}}
 \put(250,190){\line(1,0){200}}
  \put(250,190){\line(1,2){100}}
  \put(450,190){\line(-1,2){100}}
   \put(330,230){$t_2$}
\end{picture}
\end{center}
Then the derivation $D$ of $t=t'$ has three kinds of rewrite steps:
\begin{enumerate}
  \item using equation (\ref{identity2}) to the root symbol $m$ in the term;
  \item using equations (\ref{identity1}) to change something in the sequence of unary symbols over the first $\alpha$ on the left;
  \item using either kinds of equations to rewrite something in subterm $s$ or $t_2$.
\end{enumerate}
The rewrite steps of the third kind are independent of the rewrite steps of the first and second kind. Thus we can assume that we first do the rewrite steps of the first and second kind and after that the rewrite steps of the third kind. It is also not difficult to note that the rewrite steps of the first kind can be moved so that they are performed one after the other. Any two rewrite steps of the first kind done one immediately after the other do not change the term. We can eliminate all but possibly one rewrite step of the first kind from $D$ and still have a derivation of $t=t'$. Now we assume that $D$ is a derivation with at most one step of the first kind, in between the rewrite steps of the second kind, with all the rewrite steps of the third kind at the end.

Thus we have two cases depending whether there is one rewrite steps of the second kind or none. In both cases the term $t'$ is of form
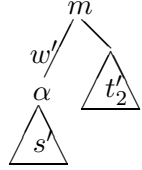
\begin{center} \xext=500 \yext=600
\begin{picture}(\xext,\yext)(\xoff,\yoff)
 \put(0,0){\line(1,0){200}}
  \put(0,0){\line(1,2){100}}
  \put(200,0){\line(-1,2){100}}
   \put(80,40){$s'$}
 \put(80,220){$\alpha$}
 \put(120,300){\line(1,2){100}}
 \put(70,350){$w'$}
   \put(200,520){$m$}
 \put(250,500){\line(1,-1){100}}
 \put(250,190){\line(1,0){200}}
  \put(250,190){\line(1,2){100}}
  \put(450,190){\line(-1,2){100}}
   \put(330,230){$t'_2$}
\end{picture}
\end{center}

If there are no rewrite step of the first kind in $D$ then  the derivation $D$ consists of three independent derivations in $T$ of the equations
\[  w=w',\;\;\;  s=s',\;\;\;  t_2=t'_2.  \]
Thus using inductive assumption we get
\[ \widehat{t} = m(u(\alpha(\widehat{s}),\widehat{t_2)} =m(u(\alpha(\widehat{s'}),\widehat{t'_2)} = \widehat{t'} \]

If there is one rewrite step of the first kind in $D$ then the derivation $D$ consists of two independent derivations in $T$ of the equations
\[  s=s',\;\;\;  t_2=t'_2.  \]
and moreover two derivation of either $w=u$ and $v=w'$ or, if $T\vdash u=v$, two derivations  $w=v$ and $u=w'$. Between the latter two derivations there is one rewrite step of the first kind. Again using inductive assumption we get
\[ \widehat{t} = m(u(\alpha(\widehat{s}),\widehat{t_2})) =m(v'(\alpha(\widehat{s'}),\widehat{t'_2})) = \widehat{t'} \]
where
\[ v'\;\;= \;\; \left\{ \begin{array}{ll}
                u   & \mbox{ if } T\vdash u=v \\
                v  & \mbox{ otherwise. }
                                    \end{array}
                            \right. \]
$\boxempty$

\noindent
\end{document}